\documentclass[12pt]{amsart}
\usepackage{amsfonts}

\newtheorem{theorem}{Theorem}

\newtheorem{example}[theorem]{Example}

\newtheorem{lemma}[theorem]{Lemma}

\newtheorem{proposition}[theorem]{Proposition}
\newtheorem{remark}[theorem]{Remark}

\begin{document}
\title[Optimality of generalized Bernstein operators]{Optimality of generalized Bernstein operators}
\author[J. M. Aldaz and H. Render]{J. M. Aldaz and H. Render}
\address{J. M. Aldaz: Departamento de Matem\'aticas, Universidad Aut\'onoma de
Madrid, Cantoblanco 28049, Madrid, Spain.}
\email{jesus.munarriz@uam.es}
\address{H. Render: School of Mathematical Sciences, University College Dublin,
Dublin 4, Ireland.}
\email{render@gmx.de}

\begin{abstract} We show that a  certain optimality property of the classical
Bernstein operator also holds, when suitably reinterpreted, for generalized
Bernstein operators on extended Chebyshev systems.
\end{abstract}

\thanks{The authors are partially supported by Grant
MTM2009-12740-C03-03 of the D.G.I. of Spain.}

\maketitle

\section{Introduction} It is well known that the classical Bernstein
operator  $B_{n}: C\left[ 0,1\right] \rightarrow
\operatorname{Span}\{1, x, \dots, x^n\},$ defined by
\begin{equation}
B_{n}f\left( x\right) =\sum_{k=0}^{n}f\left( \frac{k}{n}\right) \binom{n}{k}
x^{k}\left( 1-x\right) ^{n-k}  \label{defBP}
\end{equation}
has good shape preserving properties, a fact that explains its
usefulness in computer assisted geometric design. But the uniform
convergence of  the polynomials
$B_{n}f$ to $f$  can be
slow. Thus, it  is natural to enquire whether
one can have operators from the continuous functions to
the space of polynomials,
with good shape preserving  and fast approximation properties simultaneously.
However, in a certain sense it is impossible to do better than by using Bernstein operators.
Within some natural classes of polynomial valued operators,
Bernstein operators approximate convex functions in an optimal way.
In \cite{BeVo}, Berens and DeVore consider operators
$L_n: C\left[ 0,1\right] \rightarrow \operatorname{Span}\{1, x, \dots, x^n\},$ fixing $1$ and
$x$, and such that  $(L_nf)^{(j)}\geq 0$ if
$f^{(j)}\geq 0$ for $j=0,1,2,\cdots,n$.
They show that for every $n$, $L_n x^2 \ge B_n x^2$, and if for some
$t\in (0,1)$,
$L_n t^2 = B_n t^2$, then $L_n \equiv B_n$. Since $B_n x^2 \ge x^2$,
$L_n$ performs strictly worse than $B_n$ (unless $L_n \equiv B_n$) at least for $x^2$.
The result of Berens and DeVore \cite[p. 214]{BeVo} is extended to a wider class of operators $L_n$
in \cite[Theorem 1]{BuQue}. While the result is stated there for $x^2$,
the argument works for an arbitrary convex function $\phi$ (using
Jensen's inequality instead of H\"older's inequality). The more
general result, for $\phi$ convex, can be found in \cite[Theorem 2]{CC}. A further
extension appears in \cite[Theorem 2]{GaMi}.

 Recently, there has been an increasing interest 
concerning \emph{generalized Bernstein bases} in the framework of extended
Chebyshev spaces (in particular, spaces that contain 
functions such as $\sin x$ and $\cos x$), cf. for instance, \cite{CMP04}, \cite{CMP07}, \cite{Cost00}, \cite{CLM}, 
\cite{MPS97}, \cite{Mazu99}, \cite{Mazu05}. 
The existence of generalized Bernstein bases raises the question of the
possible existence of associated generalized Bernstein operators.
This topic is studied in \cite{AKR1},
for exponential polynomials, and in  \cite{AKR3}, \cite{AKR2}, \cite{KR},
\cite{Mazu09},
for extended Chebyshev spaces. Here we
show that  the operators considered in \cite{AKR1},
 \cite{AKR3}, \cite{AKR2}, \cite{KR}, and \cite{Mazu09}, 
 share the 
preceding optimality property of the classical Bernstein operator, under a suitable generalized notion of convexity. 

Additionally, we give a selfcontained and streamlined presentation of the generalized Bernstein operators introduced  in \cite{MoNe00}  by  S. Morigi and M. Neamtu. These
operators are  used here to exhibit  a sequence of generalized Bernstein operators converging
strongly to the identity, preserving $(1, e^x)$-convexity, and failing to
preserve standard convexity.

\section{Definitions and Optimality Results}

For each pair $n,k$ of nonnegative integers, let $\mu_{n,k}$ be
a positive measure. Denote by $p_{n,k} (x) := \binom{n}{k}
x^{k}\left( 1-x\right) ^{n-k}$ the elements of the Bernstein basis
for the space of polynomials. Let $\mathcal L_n$ be the class of positive linear operators
$L_n: C\left[ 0,1\right] \rightarrow \operatorname{Span}\{1, x, \dots, x^n\}$ fixing $1$ and
$x$ (so $L_n$ preserves affine functions) and defined by
\begin{equation}\label{defl}
L_n f (x) := \sum_{k=0}^{n} \lambda_{n,k}(f)
p_{n, k} (x), \ \ \ \ \mbox{ where } \ \ \ \  \lambda_{n,k}(f) := \int_0^1 f d \mu_{n,k}.
\end{equation}
Denote by $\delta_x$ the Dirac delta probability, which assigns 
mass 1 to the singleton $\{x\}$, and observe 
that $B_n\in \mathcal L_n$ is   obtained by taking $\mu_{n,k} = \delta_{k/n}$ in (\ref{defl}), so
$\lambda_{n,k}(f) = f(k/n)$. Recall that if $\phi$ is convex 
on $[0,1]$, then $B_n\phi \ge \phi$ on that interval.
From \cite[Theorem 1]{BuQue} and the generalization
in \cite[Theorem 2]{CC}, it is known that $L_n\phi \ge B_n\phi$ for every convex function $\phi$ and every $L_n \in \mathcal L_n$, and furthermore, $L_n\equiv B_n$ if $L_n\psi (t) =  B_n\psi (t)$
for some strictly convex
function $\psi$ and some $t\in (0,1)$. In this sense the operator $B_n$
approximates convex functions optimally within the class $\mathcal L_n$.
We emphasize  that this is not an asymptotic optimality result,
where $n$ is required to approach $\infty$. Instead, best approximation holds for each fixed
$n$.

Observe that
by uniqueness of the coefficients of the basis functions, the condition
$L_n 1 \equiv 1$ entails that each $\mu_{n,k}$ is a probability measure,
while $L_n x \equiv x$
entails that the expectation of $x$ with respect to $\mu_{n,k}$ is
$\int_0^1 x d \mu_{n,k} (x) = k/n$. If the assumption $L_n x \equiv x$ is
dropped, then the resulting operator may fail to preserve convexity,
and may approximate some convex functions better than $B_n$
(while in other cases the approximation will be worse, notably
for $x$ itself).

\begin{example}\label{ex1}{\rm Let $1 < j \le n$. In \cite[Proposition 11]{AKR3} and in \cite[Example 5.3]{Mazu09},
 a generalized Bernstein
operator $\mathcal{B}_{n,j}$  fixing $1$ and $x^{j}$,
  is defined  by
\begin{equation*}
\mathcal{B}_{n,j}f(x)=\sum_{k=0}^{n}f\left( \left[ \frac{k(k-1)\cdots (k-j+1)}{
n(n-1)\cdots (n-j+1)}\right] ^{1/j}\right) \binom{n}{k}x^{k}(1-x)^{n-k}.
\end{equation*}
In particular, if $j = 1$ this is the classical Bernstein operator.
However, for $j > 1$ and all $k = 1, 2, \dots$,
$\mathcal{B}_{n,j} x^k < B_{n} x^k$ on $(0,1)$, since the operator 
$\mathcal{B}_{n,j}$ samples the strictly increasing functions $x^k$
 at nodes lower
than those of $B_{n}$ (except  the first and the last nodes, which are
0 and 1 in both cases). It is shown in \cite[Example 17]{AKR3} that
if $1\le k\le j$, then  $\mathcal{B}_{n,j} x^k \le x^k$, while if
 $j\le k$, then $\mathcal{B}_{n,j} x^k \ge x^k$. Hence, on $(0,1)$, 
${B}_{n} x^k > \mathcal{B}_{n,j} x^k \ge  x^k$  when
 $j\le k$. Of course,  $\mathcal{B}_{n,j}$ does not preserve
 convexity, since  
 $$
 \mathcal{B}_{n,j} (1 - x) = \mathcal{B}_{n,j}(1)
 -\mathcal{B}_{n,j}( x) > 1-  B_n (x) = 1 - x
 $$ 
 whenever $x\in (0,1)$. But
 $\mathcal{B}_{n,j}$ does  preserve
$(1,x^j)$-convexity, defined below.}
\end{example}

The way to recover the results from \cite{BuQue} and \cite{CC} is thus
to deal with a generalized notion of convexity. We give next
the relevant definitions.
An \emph{extended Chebyshev
space}  $U_{n}$ of dimension $n+1$ over the interval
$\left[ a,b\right] $ is an $n+1$ dimensional subspace of
$C^{n}\left( \left[a,b\right] \right)$ such that each function $f\in U_{n}$
not vanishing identically, has at most $n$ zeros in $\left[ a,b\right]$,
counting multiplicities. It is well-known that extended Chebyshev spaces possess \emph{non-negative
Bernstein bases}, i.e. collections of non-negative functions
$p_{n,k},k=0,...,n,$ in $U_{n}$, such that each $p_{n,k}$ has a zero of order
$k$ at $a$ and a zero of order $n-k$ at $b$, for $k=0,...,n$.
Now let us select two functions $f_{0},f_{1}\in U_{n},$ such that $f_{0} >0 $ and
 $f_{1}/f_{0}$ is
strictly increasing (these functions play the role of 1 and $x$ in the
classical case). It is sometimes possible to define a generalized
Bernstein operator $\mathcal{B}_{n}:C\left[ a,b\right] \rightarrow U_{n}$
fixing $f_{0}$ and $f_{1}$, by first suitably choosing
nodes $t_{n,0},...,t_{n,n}\in \left[ a,b\right]$ such that $t_{n,0} = a$ and $t_{n,n}= b$, second, by selecting appropriate weights $\alpha _{n,0},...,\alpha _{n,n} > 0$, and finally, by
setting
\begin{equation}
\mathcal{B}_{n}\left( f\right) =\sum_{k=0}^{n}f\left( t_{n,k}\right) \alpha
_{n,k}p_{n,k}.  \label{eqBern}
\end{equation}

Criteria for the existence of such $\mathcal{B}_{n}$ can be found in
 \cite{AKR3}, \cite{AKR2}, \cite{Mazu09}, and in the specific case of exponential polynomials, in \cite{AKR1} and \cite{KR}. In this section,
we shall {\em assume} that $\mathcal{B}_{n}$ exists (fixing the given
functions), and will show
that it has an extremal property analogous to that of ${B}_{n}$.
For simplicity in the notation, we relabel the basis functions
so that they already incorporate the weights $\alpha _{n,k}$. Thus,
(\ref{eqBern}) becomes
\begin{equation}\label{eqber}
\mathcal{B}_{n}\left( f\right) =\sum_{k=0}^{n}f\left( t_{n,k}\right) p_{n,k}.
\end{equation}
 The functions $f_{0}$ and
 $f_{1}$ are always assumed to satisfy the conditions ``$f_{0} >0$" and
 ``$f_{1}/f_{0}$ is strictly increasing".

\begin{remark} {\rm Given an extended Chebyshev space $U_n$, it is
always possible to {\em select} some pair $f_0$ and $f_1$ in $U_n$ for which
a generalized Bernstein operator can be defined (cf. \cite[Corollary 8]{AKR3}).
However, for some choices of $f_0$ and $f_1$,  $\mathcal{B}_{n}$ may fail to exist. For instance, it is well known that if $0 < b < 2\pi$, then
$U_3 :=\operatorname{Span} \{1,x, \cos x, \sin x\}$ on the interval $[0,b]$
is an extended Chebyshev space. But 
 there is no generalized Bernstein operator $\mathcal{B}_{n}$ fixing $1$ and $x$  if $b=5$ (for instance) while such $\mathcal{B}_{n}$ does exist if $b = 4$, cf. \cite[Theorem 25]{AKR3}.}
\end{remark}

 \begin{remark} \label{interp}{\rm The condition on the zeros of $p_{n,k}$ entails that, like the classical Bernstein operator, the operators
 given by (\ref{eqBern}) interpolate functions at the endpoints of the interval $[a,b]$: Since $f_0$ is strictly positive, from 
$f_0 (a) = \mathcal{B}_{n}\left( f_0\right) (a) = f_0\left( a\right) \alpha
_{n, 0}p_{n,0} (a)$ we conclude that $\alpha
_{n, 0}p_{n,0} (a) = 1$. Likewise, $\alpha_{n, n}p_{n,n} (b) = 1$,
so we always have $\mathcal{B}_{n}\left( f\right) (a) = f(a)$ and $\mathcal{B}_{n}\left( f\right) (b) = f(b)$.}
\end{remark}

Convexity of $\phi$ can be defined by saying that if we interpolate
between $\phi (x)$ and $\phi(y)$ using an affine function $h$, then
on $[x,y]$ the graph of $\phi$ lies below the graph of $h$. If instead of affine functions we use $(f_0, f_1)$-affine
functions, that is, functions in
$\operatorname{Span} \{f_0, f_1\}$, we obtain the corresponding
notion of {\em $(f_0, f_1)$-convexity}.
Thus, ordinary convexity corresponds to $(1,x)$-convexity.
According
to \cite{KaSt}, p. 376, this generalized notion of convexity was introduced
 in 1926 by Hopf, and  was later
extensively developed by Popoviciu, specially in the context of Chebyshev
spaces. Strict $(f_0, f_1)$-convexity is defined analogously to strict
convexity. It is shown in \cite[Theorem 22]{AKR3} (cf. also \cite[Proposition 4.14]{Mazu09}) that a
 generalized Bernstein operator $\mathcal{B}_{n}$ fixing 
 $f_0$ and $f_1$ preserves
$(f_0, f_1)$-convexity.

We would like to stress the point  that $\mathcal{B}_{n}$ will in general {\em fail}  to preserve standard convexity (recall example \ref{ex1}).  It will be
shown in Section 3 that likewise, the standard Bernstein operator ${B}_{n}$ does not in
general preserve strengthened forms of convexity, such as $(1,e^x)$-convexity
for increasing functions.

The following characterization of $(f_0, f_1)$-convexity,
due to M. Bessen\-yei and Z. P\'ales (cf. \cite[Theorem 5]{BePa}),
helps to understand its meaning.

\begin{theorem}
\label{bepa} Let $I:= (f_1/f_0)([a,b])$. Then  $\phi\in C[a,b]$ is $(f_0, f_1)$-convex if and only
if  $(\phi/f_0) \circ (f_1/f_0)^{-1} \in C\left(I\right)$
is convex in the standard sense.
\end{theorem}

In particular, if $f_0 = 1$, the relationship between convexity
and $(1, f_1)$-convexity is given by a simple change of variables
(determined by $f_1$). For instance, returning to Example \ref{ex1},  
and to the functions considered there, we see that
if $1 \le k \le j$, then $x^k$ is $(1, x^j)$-concave, while if
 $j \le k \le n$, then $x^k$ is $(1, x^j)$-convex. 
 
Before we prove the announced optimality property of generalized Bernstein operators, we need the following lemma. It is entirely analogous to Theorem
\ref{bepa}, save that strict convexity replaces convexity throughout.
The proof is also essentially the
same, and thus  we omit it. Basically, all one needs to do is to use strict inequalities $>$ instead of $\ge$, at the appropriate places. 

\begin{lemma}
\label{bepas} Let $I:= (f_1/f_0)([a,b])$. Then  $\psi\in C[a,b]$ is strictly $(f_0, f_1)$-convex if and only
if  $(\psi/f_0) \circ (f_1/f_0)^{-1} \in C\left(I\right)$
is strictly convex.
\end{lemma}

\begin{theorem}
\label{ThmExt} Let $U_{n}\subset C^{n}\left[ a,b\right]$
be an extended Chebyshev space containing the functions $f_0$ and
$f_1$,  with $f_0 > 0$ and $f_1/f_0$ strictly increasing. Suppose
 there exists a generalized Bernstein operator
 $\mathcal{B}_{n}:C\left[ a,b\right]
\rightarrow U_{n}$ fixing   $f_{0}$ and $f_{1}$. Let the operator
 $L_{n}:C\left[ a,b\right]
\rightarrow U_{n}$, subject to $L_n f_0 = f_0$ and
$L_n f_1 = f_1$, be defined by
\begin{equation}\label{eqberns}
L_{n}\left( f\right) =\sum_{k=0}^{n}\lambda_{n,k} (f) p_{n,k},
\end{equation}
where $p_{n,k}$ has the same meaning as in (\ref{eqber}), and
$\lambda_{n,k}(f) $ is obtained from $f$ and the positive measure
$\mu_{n,k}$ via
\begin{equation}\label{mu}
\lambda_{n, k}(f) := f_0 (t_{n,k}) \int_a^b f d \mu_{n,k}.
\end{equation}
Then, for every $(f_0, f_1)$-convex function $\phi$,
we have $\phi\le \mathcal{B}_{n} \phi \le L_n \phi$. Moreover, if for
some $t\in (a,b)$ and some strictly $(f_0, f_1)$-convex function $\psi$ we have $\mathcal{B}_{n} \psi (t) = L_n \psi (t)$, then $\mathcal{B}_{n} \equiv L_n$.
\end{theorem}

\begin{proof} First we make explicit what the assumptions in the
theorem entail about
the measures $\mu_{n,k}$. Since
\begin{equation*}
L_{n} f_0 =\sum_{k=0}^{n}  {\lambda }_{n, k}( f_0 )p_{n,k} = f_0 = \mathcal{B}_{n} f_0 = \sum_{k=0}^{n}f_{0}
\left( {t_{n, k}}\right) p_{n,k},
\end{equation*}
by uniqueness of the coefficients of the basis functions, we must
have $f_{0} \left( {t_{n, k}}\right) = {\lambda }_{n, k}( f_0 )$
for all pairs $n,k$. But then $1 = \int_a^b f_0 d \mu_{n,k}$ by (\ref{mu}), so
admissible measures $\mu_{n,k}$ must be such that
$d P_{n,k}:= f_0 d \mu_{n,k}$ defines a probability $P_{n,k}$. With
this notation, the second condition $L_n (f_1) = f_1$ leads to
\begin{equation*}
{\lambda }_{n, k}( f_1 )  = f_{0}
\left( {t_{n, k}}\right) \int_a^b \frac{f_1}{f_0} d P_{n,k},
\end{equation*}
so again by uniqueness of the coefficients of the basis functions, we
have
\begin{equation*}
f_{0}
\left( {t_{n, k}}\right) \int_a^b \frac{f_1}{f_0} d P_{n,k}
=
{\lambda }_{n, k}( f_1 )  = f_{1}
\left( {t_{n, k}}\right) =  f_{0}
\left( {t_{n, k}}\right) \frac{f_{1}
\left( {t_{n, k}}\right)}{f_{0}
\left( {t_{n, k}}\right)}.
\end{equation*}
Thus,
\begin{equation}\label{aver}
\frac{f_{1}
\left( {t_{n, k}}\right)}{f_{0}
\left( {t_{n, k}}\right)}
= \int_a^b \frac{f_1}{f_0} d P_{n,k}
\end{equation}
for all pairs $n,k$.

Next, suppose   $\phi$ is $(f_0, f_1)$-convex on  $[a,b]$.
Given $h:X\to Y$ and a measure $\nu$ on $X$, the pushforward
to $Y$
 of $\nu$, using  $h$,  is denoted by $h_* \nu$; recall
that  $h_* \nu (A) := \nu (h^{-1} (A))$.
Applying the change of variables formula twice (cf., for instance,
\cite[Theorem C, p. 163]{Ha}), Theorem
\ref{bepa} above,  Jensen's
inequality, and (\ref{aver}), we obtain
\begin{equation}\label{eqjen0}
{\lambda }_{n, k}( \phi )
=
f_{0} \left( {t_{n, k}}\right) \int_a^b \phi (x) d \mu_{n,k} (x)
=
f_{0} \left( {t_{n, k}}\right) \int_a^b \frac{\phi}{f_0} (x) d P_{n,k} (x)
\end{equation}
\begin{equation}\label{eqjen1}
=
f_{0} \left( {t_{n, k}}\right) \int_{\frac{f_1}{f_0} ([a,b])} \frac{\phi}{f_0}\circ\left(\frac{f_1}{f_0} \right)^{-1} (x) d \left(\frac{f_1}{f_0} \right)_* P_{n,k} (x)
\end{equation}
\begin{equation}\label{eqjen2}
\ge
f_{0} \left( {t_{n, k}}\right) \frac{\phi}{f_0}\circ\left(\frac{f_1}{f_0} \right)^{-1} \left(\int_{\frac{f_1}{f_0} ([a,b])} x
d \left(\frac{f_1}{f_0} \right)_* P_{n,k} (x)\right)
\end{equation}
\begin{equation}\label{eqjen3}
=
f_{0} \left( {t_{n, k}}\right) \frac{\phi}{f_0}\circ\left(\frac{f_1}{f_0} \right)^{-1} \left(\int_{[a,b]} \left(\frac{f_1}{f_0} \right) (x)
d  P_{n,k} (x)\right)
\end{equation}
\begin{equation}\label{eqjen4}
=
f_{0} \left( {t_{n, k}}\right) \frac{\phi}{f_0}\circ\left(\frac{f_1}{f_0} \right)^{-1} \left(\frac{f_1}{f_0} (t_{n,k})\right) = \phi (t_{n,k}).
\end{equation}
Since the basis functions are non-negative,
we obtain $\mathcal{B}_{n} \phi \le L_n \phi$ just by adding up.
And the fact that $\phi\le \mathcal{B}_{n} \phi$
for all
$(f_0, f_1)$-convex functions $\phi$ is proven in \cite[Theorem 15]{AKR3}.

Assume next that for some $t\in (a,b)$ and some strictly $(f_0, f_1)$-convex
function $\psi$ we have $\mathcal{B}_{n} \psi (t) = L_n \psi (t)$.
We must show that  $P_{n,k} = \delta_{t_{n,k}}$
for each pair $n, k$. Set  $\phi =\psi$
in (\ref{eqjen0})-(\ref{eqjen4}).
Since
\begin{equation*}
\sum_{k=0}^{n}  {\lambda }_{n, k} ( \psi )p_{n,k} (t)
=
L_{n} \psi (t)
=
\mathcal{B}_{n} \psi (t)
 =
 \sum_{k=0}^{n}\psi
\left( {t_{n, k}}\right) p_{n,k} (t),
\end{equation*}
and for all $n,k$, we have ${\lambda }_{n, k} ( \psi ) \ge \psi
\left( {t_{n, k}}\right)$ (by (\ref{eqjen0})-(\ref{eqjen4})) from ${p}_{n, k} > 0$ on
$(a,b)$ we conclude that ${\lambda }_{n, k} ( \psi ) = \psi
\left( {t_{n, k}}\right)$ for all $n,k$. Thus, with $\phi =\psi$ we have equality in (\ref{eqjen2}).
By Lemma \ref{bepas}, the function
 $\frac{\psi}{f_0}\circ\left(\frac{f_1}{f_0} \right)^{-1}$ is strictly convex.
So by the equality case in Jensen's inequality 
we conclude that the function $x$ is constant a.e. with respect to
$\left(\frac{f_1}{f_0} \right)_* P_{n,k}$. This entails that the latter
measure is a Dirac delta, and since  ${f_1}/{f_0}$ is strictly increasing,
 $P_{n,k}$ must also be a Dirac delta. It now follows from (\ref{aver})
that $P_{n,k} = \delta_{t_{n,k}}$.
\end{proof}

\section{On the generalized Bernstein Operator of Morigi and Neamtu}

There are known examples of chains of extended Chebyshev spaces,
for which
 their associated generalized Bernstein operators 
  exist and converge strongly 
to the identity. An instance is given by the 
operators $\mathcal{B}_{n,j}$ described in Example \ref{ex1}  (cf.  \cite[Proposition 11]{AKR3}  or \cite[Theorem 6.1]{Mazu09}
for the convergence assertion). As was observed in Example \ref{ex1}, the operators 
$\mathcal{B}_{n,j}$ preserve $(1,x^j)$-convexity, and fail to preserve ordinary convexity when $j > 1$.

Nevertheless, these operators are somewhat
degenerate in that the sequences of nodes are not strictly increasing,  
due to the fact that for $j>1$, $f_1(x) = x^j$ has a zero of order at least 
two
at zero, and in particular,  $f^\prime_1(0) = 0$. Thus, $f^\prime_1$ is not strictly positive on $[0,1]$. We shall
use the generalized Bernstein operators introduced by Morigi and Neamtu
in \cite{MoNe00} to present a sequence of operators that i) converges
to the identity, ii) the operators preserve $(1, e^x)$-convexity, iii) they do not preserve ordinary convexity, iv) they are defined via strictly increasing
sequences of nodes (in fact, the nodes are equidistributed), and v)
they fix $f_0 = 1$ and   $f_1(x)/f_0(x) = f_1(x) = e^x$. Obviously we have $f_1^\prime > 0$  everywhere.

Our presentation here of the  Morigi and Neamtu operators is selfcontained and
simplifies some of the arguments from \cite{MoNe00}. At the beginning we
will consider a setting more general than we really need (from the
shape preservation perspective) allowing complex values. We start by defining these operators.

Given two different complex numbers $\mu_{0},\mu_{1}$, set
$\omega_{n}:=\frac{1}{n}\left(  \mu_{1}-\mu_{0}\right)  \neq0$, and for 
$j=0,...,n$, let 
\begin{equation}\label{exponen}
\lambda_{j}:=\mu_{0}+j\omega_{n}.
\end{equation}
Thus,  $\lambda_{0}=\mu_{0}$ and $\lambda_{n}=\mu_{1}.$ Define
\begin{equation}
\varphi_{n}\left(  x\right)  :=\frac{e^{\omega_{n}\left(  x-b\right)
}-e^{\omega_{n}\left(  a-b\right)  }}{1-e^{\omega_{n}\left(  a-b\right)  }
}\text{ and }\psi_{n}\left(  x\right)  
:=
1 - \varphi_{n}\left(  x\right)
=
\frac{1-e^{\omega_{n}\left(
x-b\right)  }}{1-e^{\omega_{n}\left(  a-b\right)  }}.\label{eqphin}
\end{equation}
Then $\left(  \varphi_{n}\right)  ^{k}$ is a linear combination
of the exponential functions $e^{\omega_{n}jx}$, where $j=0,...,k,$ and $\left(
\psi_{n}\right)^{n-k}$ is a linear combination of the functions $e^{\omega_{n}jx}$, where 
$j=0,...,n-k,$. Therefore, 
\begin{equation}\label{basisMN}
p_{n,k}\left(  x\right)  := e^{-\lambda_{0}\frac{k}{n}\left(  b-a\right)
}\cdot e^{\lambda_{0}\left(  x-a\right)  }\binom
{n}{k} \left(  \varphi_{n}\left(  x\right)
\right)  ^{k} \left(  \psi_{n}\left(  x\right)
\right)  ^{n-k}
\end{equation}
is a linear combination of the exponential functions $e^{\lambda_{j}x}$, $j=0,...,n$.
It follows that $p_{n,k}$ is
contained in the vector space
\begin{equation}
U_{n}:=\left\{  f\in C\left(  \mathbb{R},\mathbb{C}\right)  :\left(  \frac
{d}{dx}-\lambda_{0}\right)  \cdots\left(  \frac{d}{dx}-\lambda_{n}\right)
f=0\right\}.
\end{equation}
Since $\varphi_{n}$ has a simple zero at $x=a$ and $\psi_{n}$ a simple zero at
$b$, $p_{n,k}$ has a zero order $k$ at $a$ and a zero of
order $n-k$ at $b.$ 

It is well known that if $\lambda_{0},...,\lambda_{n}$ are real, then $U_{n}^{r}:=\left\{  f\in
U_{n}:f\text{ is real-valued}\right\}  $ is an extended Chebyshev space over
any bounded closed interval $\left[  a,b\right]$
(we mention that when dealing with
 complex values $\lambda_{0},...,\lambda
_{n}$, restrictions must be imposed on the length of $[a,b]$ in order to ensure that $U_{n}$ is an extended Chebyshev space, cf.
for instance
\cite[Theorem 19]{AKR1}). Furthermore, $\{p_{n,k}: k=0, \dots, n\}$ is a
non-negative Bernstein basis. 

Next we specify nodes and weights. Let  $t_{n,k}:=a+\frac{k}{n}\left(  b-a\right)  $ for $k=0,...,n$, and let
\begin{equation}\label{bopmone}
\mathcal{B}_{n}\left(  f\right)  \left(  x\right)  :=\sum_{k=0}^{n}
f\left(  t_{n,k}\right) p_{n,k}\left(  x\right).
\end{equation}

\begin{proposition} The operator $\mathcal{B}_{n}$ fixes $e^{\mu_{0}x}$
and $e^{\mu_{1}x}$.
\end{proposition}

\begin{proof} Recall from (\ref{exponen}) that $\mu_0 = \lambda_0$ and $\mu_1 =
\lambda_n$. Using (\ref{bopmone}), (\ref{basisMN}), and (\ref{eqphin}), we obtain
\[
\mathcal{B}_{n}\left(
e^{\lambda_{0}x}\right)  =
e^{\lambda_{0}x} \sum_{k=0}^{n}\binom{n}{k}\left(
\varphi_{n}\left(  x\right)  \right)  ^{k}\left(  1-\varphi_{n}\left(
x\right)  \right)  ^{n-k}
=e^{\lambda_{0}x} \left(\varphi_{n}\left(  x\right) +  1-\varphi_{n}\left(
x\right)  \right)  ^{n}
=e^{\lambda_{0}x}.
\]
Furthermore,
\begin{align*}
\mathcal{B}_{n}\left(  e^{\lambda_{n}x}\right)  &
=
\sum_{k=0}^{n}
e^{\lambda_{n}a}e^{\left(  \lambda_{n}-\lambda_{0}\right)  
\frac{k}{n}\left(  b-a\right)  }e^{\lambda_{0}\left(  x-a\right)  }
\binom{n}{k}\left(  \varphi
_{n}\left(  x\right)  \right)  ^{k}\left(  1-\varphi_{n}\left(  x\right)
\right)  ^{n-k}\\
& =
e^{\lambda_{n}a}e^{\lambda_{0}\left(  x-a\right)  }\sum_{k=0}^{n}\binom
{n}{k}\left(  e^{\omega_{n}\left(  b-a\right)  }\varphi_{n}\left(  x\right)
\right)  ^{k}\left(  1-\varphi_{n}\left(  x\right)  \right)  ^{n-k}\\
& =
e^{\lambda_{n}a}e^{\lambda_{0}\left(  x-a\right)  }\left(  e^{\omega
_{n}\left(  b-a\right)  }\varphi_{n} \left(  x\right)
 + 1-\varphi_{n}\left(  x\right)  \right)
^{n}\\
& =
e^{\lambda_{n}a}e^{\lambda_{0}\left(  x-a\right)  }\left(  \left(
e^{\omega_{n}\left(  b-a\right)  }-1\right)  \varphi_{n}\left(  x\right) +1\right)  ^{n}.
\end{align*}
But now (\ref{eqphin}) implies that
\begin{align*}
\left(  e^{\omega_{n}\left(  b-a\right)  }-1\right)  \varphi_{n} \left(  x\right) +1  &
=e^{\omega_{n}\left(  b-a\right)  }\left(  1-e^{\omega\left(  a-b\right)
}\right)  \varphi_{n} \left(  x\right) +1\\
& =e^{\omega_{n}\left(  b-a\right)  }\left(  e^{\omega_{n}\left(  x-b\right)
}-e^{\omega_{n}\left(  a-b\right)  }\right)  +1=e^{\omega_{n}\left(
x-a\right)  },
\end{align*}
so from the preceding equalities we obtain
\[
\mathcal{B}_{n}\left(
e^{\lambda_{n}x}\right)  =e^{\lambda_{n}a}e^{\lambda_{0}\left(  x-a\right)
}e^{n\omega_{n}\left(  x-a\right)  }=e^{\lambda_{n}x}.
\]
\end{proof}

For simplicity, next we take $\left[  a,b\right]  $ to be 
$[0,1]$, $\mu_0 = 0$, and $\mu_1 = 1$. Then the
Morigi and Neamtu operator reduces to
\begin{equation}\label{bopmonesp}
\mathcal{B}_{n} f  \left(  x\right)  :=\sum_{k=0}^{n}
f\left(  \frac{k}{n}\right) \binom
{n}{k} \left(  \varphi_{n}\left(  x\right)
\right)  ^{k}\left(  1-\varphi_{n}\left(  x\right)  \right)  ^{n-k},
\end{equation}
where
by specialization of (\ref{eqphin}) we have
\[
\varphi_{n}\left(  x\right)  
=\frac{e^{x/n}-1}{e^{1/n}-1}.
\]
Observe that since $e^{x/n}$ is convex and increasing, so is $\varphi_{n}$.
Furthermore, $\varphi_{n}\left(  0\right)  =0$ and $\varphi_{n}\left(  1\right)
=1$. Recalling that $B_n$ denotes the classical Bernstein operator,
given by  (\ref{defBP}), we see that
\begin{equation}
\mathcal{B}_{n} f \left(  x\right)  = B_{n} f
\left(  \varphi_{n}\left(  x\right)  \right).   \label{eqnice}
\end{equation}
From the preceding expression it is immediate that $\mathcal{B}_{n}$ 
does not in general preserve the convexity of a (decreasing) convex
function. Let $h(x) = -x$. Since $B_n h = h$, we have 
\[
\mathcal{B}_{n} h  \left(  x\right)  =B_{n} h
\left(  \varphi_{n}\left(  x\right)  \right)  = - \varphi_{n}\left(  x\right).
\]
But $\varphi_{n}$ is strictly convex, so the image under $\mathcal{B}_{n}$
of the affine (hence convex) function $h$ is the strictly concave
function $- \varphi_{n}$. This shows that the following assertion, 
contained in \cite[Remark 4.15 (1)]{Mazu09}, is not correct: If $f_0 = 1$,
$f_1$ is convex, 
and $f_1^\prime > 0$ everywhere on $[a,b]$,  
then the image of a decreasing convex function under $\mathcal{B}_{n}$
is convex (this is claimed to follow from \cite[Proposition 4.14]{Mazu09}, 
but it does not; we note that Proposition 4.14 
from \cite{Mazu09} is correct, cf. 
\cite[Theorem 22]{AKR3}  together with Theorem 4 above for a stronger 
result). 

It is shown in \cite[Theorem 3.3]{MoNe00} that the operators 
$\mathcal{B}_{n}$ converge strongly to
the identity (more general or related results can be found in \cite[Theorem 23]{AKR1} and
\cite[Theorem 6.1]{Mazu09}). However,  for the special case we are considering, it is easy to give
a direct argument, so in order to make this section as selfcontained as
possible, we do this next.

Let $g(x) := x$. Since $\mathcal{B}_{n}$ fixes $1$ and $e^x$, by Korovkin's Theorem it suffices to prove  that $\mathcal{B}_{n} g$
converges uniformly to $g$ on $[0,1]$. By (\ref{eqnice}), 
$\mathcal{B}_{n} g    = 
 \varphi_{n}$, so
it is enough to show that $\varphi_{n}$ converges uniformly to $g$. Since $\varphi_{n}$ is convex and increasing, its derivative achieves
its maximum value over $[0,1]$ at $1$. Now let $g_n$ be the affine function
with slope $\varphi_{n}^\prime (1)$ passing thorough the point $(1,1)$.
By convexity, $g_n \le \varphi_{n} \le g$, so to obtain the uniform
convergence it suffices to show that $\lim_{n\to\infty }\varphi_{n}^\prime (1) = 1$. But this follows immediately by writing explicitly $\varphi_{n}^\prime (1)$, and then using L'Hospital rule, or the Taylor expansion for $e^{1/n}$.

We have seen that $\mathcal{B}_{n}$ does not preserve convexity. 
It is easy to check that $B_n$ does not preserve $(1, e^x)$-convexity
either. Actually, both of these facts can be derived from the slightly more general result
given next. 

\begin{proposition} \label{noco} Let $f_0 =1$ and $f_1$ be
functions in an Extended Chebyshev Space $E$ over the
interval $[a,b]$, where $f_1$ is
assumed to be strictly increasing and strictly convex. 
Suppose there exists a generalized Bernstein operator
$\mathcal{B}_{n}: C[a,b]\to E$  fixing
$f_0$ and $f_1$. Then the classical Bernstein operator ${B}_{n}$ does not preserve  $(1, f_1)$-convexity,
and the generalized operator $\mathcal{B}_{n}$ does not preserve convexity.
\end{proposition}

\begin{proof} For the first assertion, observe that while $f_1$ is
trivially  $(1, f_1)$-convex,  $B_n f_1$ is not.  To see why, note
that
since $B_n f_1 (a) = f_1 (a)$ and $B_n f_1 (b) = f_1 (b)$,  the unique
function in $\operatorname{Span}\{1, f_1\}$ interpolating $B_n f_1$
at $a$ and $b$ is $f_1$ itself. If $B_n f_1$ were $(1, f_1)$-convex,
we would have $B_n f_1 \le f_1$ on $[a,b]$. But $f_1$ is 
strictly convex, so
by Jensen's inequality (including the equality condition)
we have $B_n f_1 (x) > f_1(x)$ for all $x\in (a,b)$.

An analogous argument  yields the second assertion. Set $h(x) = x$, so 
 both  $h$ and $-h$
are convex. Towards a contradiction, suppose  that  $\mathcal{B}_{n}$ preserves convexity.
Then $\mathcal{B}_{n} h \le h$ and 
$\mathcal{B}_{n}(-  h) \le -h$, so
 $\mathcal{B}_{n} h = h$. Fix $x\in (a,b)$, and note that 
 $f\mapsto \mathcal{B}_{n} f(x)$ is a positive linear functional
 on $C[a,b]$, defined by a probability measure 
 (since $\mathcal{B}_{n} 1 (x) = 1$). Thus we can once more apply Jensen's inequality, in this case 
 to the strictly concave function $f_1^{-1}$, to conclude that 
 $$x = h(x) = \mathcal{B}_{n} h (x) = \mathcal{B}_{n} (f_1^{-1}\circ f_1) (x)
 >
f_1^{-1} \left( \mathcal{B}_{n} f_1  (x)\right) = f_1^{-1}\left( f_1 (x)\right)
= x.
 $$
 Alternatively, we can reach the same conclusion  without using Jensen's
 inequality, by observing that if in addition to fixing $1$ and $f_1$, we have $\mathcal{B}_{n} h = h$, then $\mathcal{B}_{n}$ 
 fixes three linearly independent functions, and thus it must be
 the identity operator (by  \cite[Proposition 3.7]{Mazu09}). So again we
 have a contradiction.
\end{proof}

We finish with an easy Proposition, which sheds some light on the relationship between the different notions of convexity considered above.

\begin{proposition} \label{noco} Let $f_1$ be increasing and strictly convex
on $[a,b]$. For increasing functions on $[a,b]$,
the condition of
 $(1, f_1)$-convexity is strictly stronger than convexity.
Over the  decreasing functions on $[a,b]$,
convexity is strictly stronger than
 $(1, f_1)$-convexity.
\end{proposition}

\begin{proof} Let $f$ be increasing and 
$(1, f_1)$-convex on $[a,b]$. Pick $x,y\in [a,b]$ with $x < y$. Since $f(y)\ge f(x)$, the
unique interpolant  $\psi$ of $f$ at $x$ and $y$ such that $\psi\in\operatorname{Span}\{1, f_1\}$
must be of the form $\psi = c_0 + c_1 f_1$, with $c_1 \ge 0$.
Thus $\psi$ is convex, so the unique interpolant $\phi \in\operatorname{Span}\{1, t\}$ of $\psi$ (or $f$) at $x$ and $y$
satisfies $\phi (t) \ge \psi (t) \ge f(t)$ for all $t\in [x,y]$,
whence $f$ is convex. Also, $h(t) := t$ is increasing and
convex but not $(1, f_1)$-convex, since $h(f^{-1}_1) = f^{-1}_1$
is strictly concave, cf. Lemma \ref{bepas}. 

The assertion about decreasing functions is proved by the same type
of argument, so we omit it.
\end{proof}

\end{document}